\documentclass{article}

\usepackage{arxiv}

\usepackage[utf8]{inputenc} 
\usepackage[T1]{fontenc}    
\usepackage{hyperref}       
\usepackage{url}            
\usepackage{booktabs}       
\usepackage{amsfonts}       
\usepackage{nicefrac}       
\usepackage{microtype}      
\usepackage{lipsum}
\usepackage{graphicx}
\usepackage{amsmath}
\graphicspath{ {./images/} }

\title{SuperMeshing: A Novel Method for Boosting the Mesh Density in Numerical Computation within 2D Domain}

\author{
    \textbf{Handing Xu$^{1,3}$} \quad
    \textbf{Zhenguo Nie$^{1,2,}$}\thanks{Address all correspondences to these authors ( zhenguonie@tsinghua.edu.cn).} \quad
    \textbf{Qingfeng Xu$^{1,4}$} \quad 
    \textbf{Xinjun Liu$^{1,2,*}$}  \\ 
    
    ${^1}$The State Key Laboratory of Tribology \& \\Tsinghua University (DME)-Siemens Joint Research Center for Advanced Robotics,\\ Department of Mechanical Engineering (DME), Tsinghua University, Beijing 100084, China\\
    ${^2}$Beijing Key Lab of Precision/Ultra-precision Manufacturing Equipments and Control,\\ Tsinghua University, Beijing 100084, China\\
    ${^3}$School of Mechatronical Engineering, Beijing Institute of Technology, Beijing 100084, China\\
    ${^4}$School of Computing and Information Systems. Melbourne School of Engineering,\\ University of Melbourne, Melbourne 3000, Australia\\
}
\begin{document}
\maketitle
\begin{abstract}
Due to the limit of mesh density, the improvement of the spatial precision of numerical computation always leads to a decrease in computing efficiency. Aiming at this inability of numerical computation, we propose a novel method for boosting the mesh density in numerical computation within 2D domain. 
Based on the low mesh-density stress field in 2D plane strain problem computed by the finite element method, this method utilizes a deep neural network named SuperMeshingNet to learn the non-linear mapping from low mesh-density to high mesh-density stress field, and realizes the improvement of numerical computation accuracy and efficiency simultaneously. We adopt residual dense blocks to our mesh-density boost model – SuperMeshingNet for extracting abundant local features and enhancing the prediction capacity of the model.
Experimental results show that the SuperMeshingNet proposed in this work can effectively boost the spatial resolution of the stress field under the multiple scaling factors: $2\times, 4\times, 8\times$. Compared to the results of finite element method, the predicted stress field error of SuperMeshingNet is only $0.54\%$, which is within the acceptable range of stress field estimation, and the SuperMeshingNet predicts the maximum stress value also without significant accuracy loss.
We publicly share our work with full detail of implementation at \url{https://github.com/zhenguonie/2021_SuperMeshing_2D_Plane_Strain}.
\end{abstract}


\section{Introduction}
Numerical computation, especially the finite element analysis (FEA) is widely used in various fields of scientific research\cite{FEA_1, FEA_2}, especially in mechanical fields such as structural performance evaluation\cite{FEA_3} and topology optimization\cite{FEA_4}. Essentially, FEA is a kind of computer technology, which discretizes the continuous physical system and solves this discretization model. 
Because it determines the fineness of the continuous system discretization, mesh density is an important factor in FEA. At the same time, the computing cost of solving the discretized linear equations is linearly proportional to the mesh density. This relation results in the impossibility of improving the spatial resolution and computing efficiency of FEA simultaneously.

In view of the low computational efficiency of FEA under high mesh density, two kinds of solutions were proposed in recent years: adaptive mesh refinement (AMR), machine-learning based physical field prediction.
AMR uses coarser mesh in unimportant area, and finer mesh in the interesting area, so as to obtain a better solution with less mesh. Berger et al.\cite{AMR_1} first proposed an algorithm based on recursive iteration to achieve the local mesh density increase. After that, more research\cite{AMR_2, AMR_3, AMR_4, AMR_5} was performed to propose various enhanced AMR algorithms on this basis, and plenty of open source libraries were established\cite{AMR_6, AMR_7}. However, the AMR algorithms still require multiple iterations to gradually decide which part needs a finer mesh and obtain higher numerical accuracy. The increase in spatial resolution is still accompanied by the decrease in computing efficiency. The AMR algorithm cannot fundamentally solve the contradiction between the spatial resolution and the computing efficiency of FEA.

Machine-learning based methods can predict the physical field in one-shot getting rid of the iterative process, using the data-driven framework to map from input physical conditions to output physical fields. Liang et al.\cite{DL_FEA_1} first proposed a groundbreaking stress field estimation method based on deep learning. Their work concentrated on the stress field of the human thoracic aorta. The trained model can output the stress field of the thoracic aorta within one second after inputting the shape of the aorta. This method has significantly improved the efficiency, but only suitable for the estimation of the stress field of the thoracic aorta. Bhatnagar et al.\cite{DL_FEA_4} proposed an approximation model based on a convolutional neural network (CNN) for flow field predictions. Compared with the finite element method (FEM), the efficiency of this model is improved by 1-2 orders of magnitude. These progress proves the potential of machine learning in solving the contradiction between spatial resolution and computing efficiency. However, numerous training data are needed for these models and limit the prediction on prior unseen cases. 

\begin{figure}[!ht]
    \centering
    \includegraphics[width=0.6\linewidth]{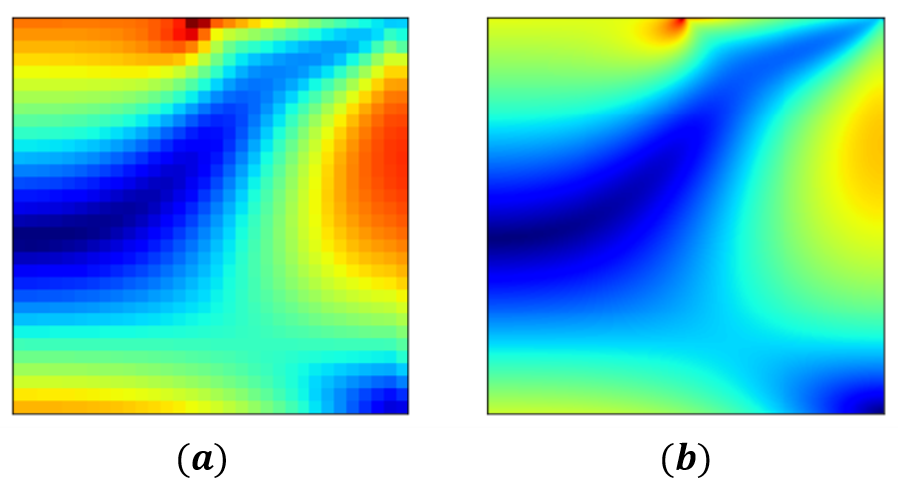}
    \caption{Low mesh-density $(32\times32)$ stress field $(a)$ and High mesh-density $(256\times256)$ stress field $(b)$.}
    \label{fig:low_high_field}
\end{figure}

In our work, we proposed a novel method for boosting the mesh density in physical field prediction within 2D domain. Our method considers the mesh density in the numerical computing field to the image resolution in the image super-resolution field. This image-based problem representation allows the analysis of arbitrary 2D planar structures. We propose the SuperMeshingNet, which focuses on the 2D plane strain problem. 
SuperMeshingNet inputs low mesh-density stress field computed by FEM, takes high mesh-density physical field as the target, and outputs the generated high mesh-density physical field (Figure \ref{fig:low_high_field} shows a pair of low mesh-density and high mesh-density stress field). Predicting the stress field by SuperMeshingNet, less computing resources are needed on high mesh density compared to FEM, thereby fairly avoiding the conflict between the spatial resolution and the computing efficiency of FEA. In our test, we generate a dataset containing 3240 different boundary and load conditions. Under the mesh-density scaling factor is $8$, the average relative stress error is only $0.54\%$, and the average computing speed is increased by $13.36$ times than FEM.

\section{Related Work}
Our review focuses on studies that highlight physical field prediction based on deep learning, image super resolution, and two neural network architectures closely related to our work.

\subsection{Physical Field Prediction Based on Deep Learning}
Lots of research\cite{DL_FEA_1, DL_FEA_2} have proved that machine learning has excellent application potential and prospects in the field of physical field prediction. Ladický et al.\cite{DL_FEA_3} proposed a regression forest-based fluid simulation method for the prediction of the flow field. Compared with FEM, the computational efficiency of this model is improved by 1-3 orders of magnitude, which has a good alternative potential of the traditional method. Bhatnagar et al.\cite{DL_FEA_4} proposed a flow field prediction model based on CNN. In the research of Nie et al., CNN and GAN have been used for structural analysis\cite{DL_FEA_5, DL_FEA_9} and topology optimization\cite{DL_FEA_8}. Besides, a physical field predicts method based on GNN had been proposed in the research of Pfaff et al.\cite{DL_FEA_2}, this model can widely apply to the dynamic simulation of various physical systems. Compared with FEM, all of the deep learning based studies have greatly improved the computing efficiency. 

Compared with machine learning completely replacing finite element analysis, another idea of combining machine learning with numerical computing has attracted more and more attention in recent years. Xie et al.\cite{DL_FEA_6} proposed a super-resolution model of flow field based on a generative adversarial network (GAN) - tempoGAN, the model can output clearer four-dimensional ($x, y, z$ and $t$) flow field by inputting the flow field obtained by FEM. Belbute et al.\cite{DL_FEA_7} combined GNN and partial differential equation (PDE) solution, proposed a flow field prediction method, calculate low mesh-density flow field by PDE, then input it with high density mesh into GNN, output the high mesh-density flow field. Although both methods obtain flow field results with high accuracy and high computing efficiency, there are still issues with versatility.

\subsection{Image Super Resolution}
Image super resolution solves the problem of HD reconstruction of low resolution images, which requires to recover natural and clear textures from low resolution images, and finally obtains high resolution images. In recent years, supervised neural networks have become the most advanced technology in the field of image super resolution. The SRCNN proposed by Dong et al.\cite{SR_1} applied CNN to image super resolution field for the first time. Shi et al.\cite{SR_2} proposed a sub-pixel convolutional layer instead of pre-upsampling in the ESPCN model, this architecture not only achieves better SR effect, but also has a great advantage in speed. Lim et al.\cite{SR_3} proved the MAE loss function is better than MSE in the field of image super resolution in their work. The VDSR model proposed by Kim et al.\cite{SR_4} and the SRDenseNet model proposed by Tong et al.\cite{SR_5} first introduced residual learning and dense block to image super resolution separately, both of these work bring great contributions to the image super-resolution field. On the basis of these models, researchers in this field have developed a variety of super resolution models\cite{SR_6, SR_7, SR_8, SR_9}.

Compare the mesh density increase of stress field and the super resolution of image, it's easy to find out that they are very similar. Therefore, using image processing methods to process the physical field is of great significance both in theory and practice.

\subsection{Residual Learning and Residual Dense Block}
\begin{figure}[!ht]
    \centering
    \includegraphics[width=0.4\linewidth]{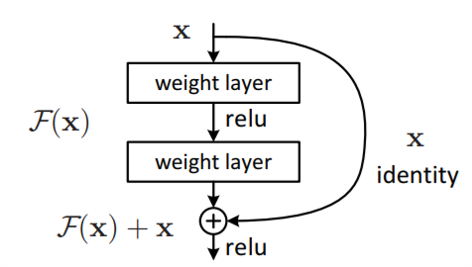}
    \caption{Residual learning architecture.\cite{Resnet}}
    \label{fig:residual_learning}
\end{figure}

The increase in the depth of the neural network will not only improve the non-linear expressiveness, it will also bring about problems such as network degradation. He et al.\cite{Resnet} proposed residual learning in order to solve this problem, the architecture is shown in Figure \ref{fig:residual_learning}. For a stacked convolutional layer model, the input is $x$, assume the feature to be learned is $H(x)$, then we let the model learn its residual $F(x)=H(x)-x$ rather than itself $H(x)$. Intuitively, because the residuals to be learned are generally smaller than original feature, there is less content to learn, result in the learning difficulty decrease.

\begin{figure}[!ht]
    \centering
    \includegraphics[width=0.75\linewidth]{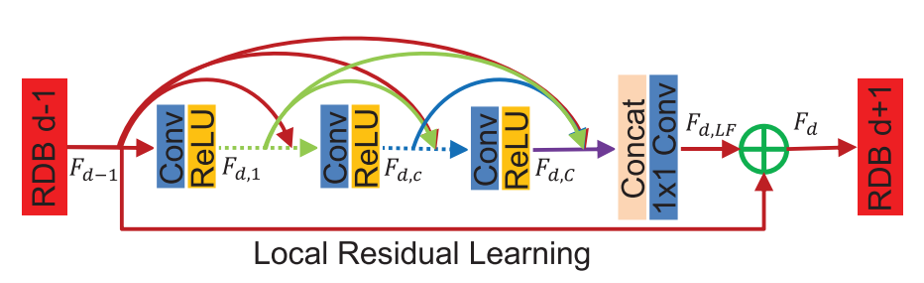}
    \caption{Residual dense block (RDB) architecture.\cite{SR_6}}
    \label{fig:residual_dense_block}
\end{figure}

Zhang et al.\cite{SR_6} proposed residual dense block (RDB) in RDN, this architecture integrates dense connection layer, local feature fusion and local residual learning. As shown in Figure \ref{fig:residual_dense_block} and Equation (\ref{equ:rdn}), it synthetically uses the features of $F_{d-1}, F_{d,c}$ and forms a continuous memory mechanism, thus can more effectively use the feature information of the image and enhance the expressive ability of the model.

\begin{equation}
    \begin{aligned}
        &F_d =  F_{d-1} + F_{d, LF} \\
        &F_{d, LF} =  H^d_{LFF}([F_{d-1}, F_{d, 1}, \cdots, F_{d, c}, \cdots, F_{d, C}]) \\
        &F_{d, c} = \sigma(W_{d, c}[F_{d-1}, F_{d, 1}, \cdots, F_{d, c-1}])
    \end{aligned}\label{equ:rdn}
\end{equation}
\noindent
where $H^d_{LFF}$ denotes the function of the $1\times1$ convolutional layer, $\sigma$ denotes the ReLU activation function and $W_{d, c}$ denotes the weights of the $c-th$ convolutional layer.

\subsection{Sub-pixel Convolution Layer}
\begin{figure}[!ht]
    \centering
    \includegraphics[width=1\linewidth]{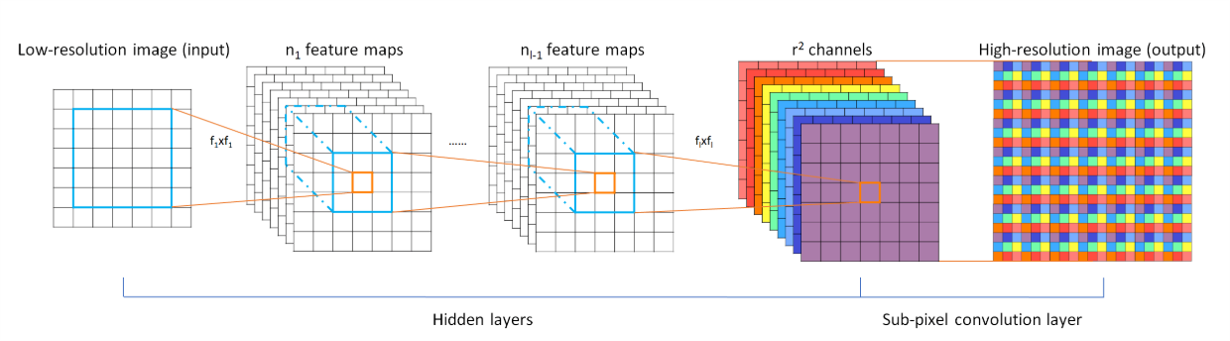}
    \caption{The sturcture of ESPCN, the last two layers make up the sub-pixel convolution layer.\cite{SR_2}}
    \label{fig:sub_pixel_conv_layer}
\end{figure}

The sub-pixel convolution layer is proposed for upsampling in ESPCN by Shi et al.\cite{SR_2} The layer architecture is shown in Figure \ref{fig:sub_pixel_conv_layer}. Assuming that the size of the low resolution image is $x\times y\times c$, the scale is $r$, the size of feature map $T$ learned in the penultimate layer of the model is $x\times y\times r^2c$. In sub-pixel convolution layer, rearranges the $r^2c$ channels to the size of $rx\times ry\times c$, the rearrangement rules are shown in Equation (\ref{equ:sub-pixel conv layer}): 

\begin{equation}
    PS(T)_{x,y,c}=T_{\lfloor{x/r}\rfloor, \lfloor{y/r}\rfloor, c\cdot r\cdot mod(y,r)+c\cdot mod(x,r)}\label{equ:sub-pixel conv layer}
\end{equation}

Though this layer architecture is called sub-pixel convolutional layer, it actually does not perform the convolution operation. It's an efficient and parameter-free pixel rearrangement method. Therefore, it greatly improves the study speed of the entire model. On the other hand, because of this architecture, the model is upsampled at the end, so that the model can retain more texture information in the low resolution space and improve the reconstruction effect.

\section{Technical Approach}
This work proposes a method for boosting the mesh density in physical field prediction within 2D domain. We calculate the stress field of 2D plane structure with linear isotropic elastic material. The approach takes the structure geometry, external loads, displacement boundary conditions and low density mesh as the input of FEM, and the low mesh-density von-mises stress field obtained is used as the input of the neural network to output the high mesh-density von-mises stress field. Based on this method, we propose a neural network architecture —— SuperMeshingNet.

\subsection{Problem Description}
\begin{figure}[!ht]
    \centering
    \includegraphics[width=0.35\linewidth]{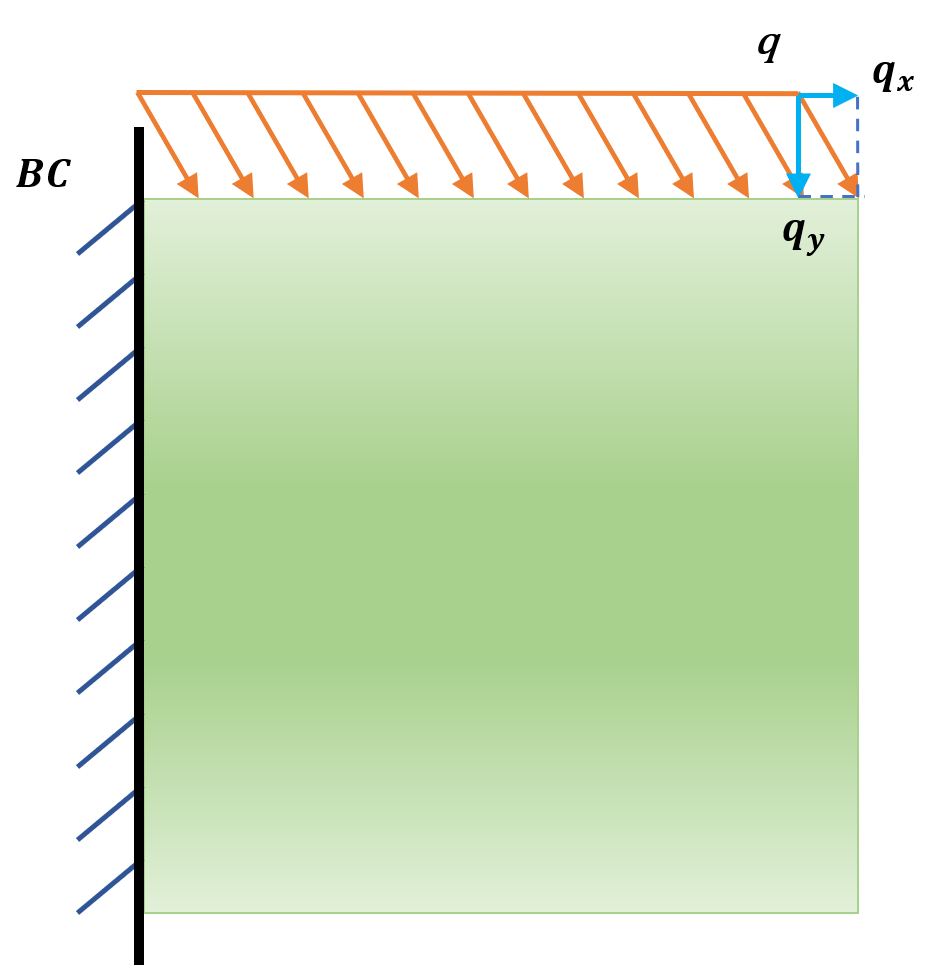}
    \caption{A sample of a plane strain problem with linear isotropic material. $BC$ denotes the boundary condition, $q$ denotes the external force.}
    \label{fig:plane_strain_problem}
\end{figure}

Consider the structure shown in Figure \ref{fig:plane_strain_problem} is composed of a homogeneous and isotropic linear elastic material. We treat this structure as a plane strain problem. The load is applied to the upper boundary of the structure, and different types of displacement boundary conditions are applied to other boundaries. The evenly distributed external force $q$ is applied in both the horizontal and vertical directions. The magnitudes of the two components denoted by $q_x$ and $q_y$ vary separately within a prescribed window of magnitude. The displacement boundary conditions includes completely fixed, horizontal constraint, vertical constraint and so on. The material properties are kept unchanged for all samples.

The constitutive equation of the plane strain problem consist of the geometric equation (\ref{equ:geometric_equation}), equilibrium equation (\ref{equ:equilibrium_equation}) and physical equation, the generalized Hooke law (\ref{equ:physical_equation}): 

\begin{equation}
    \varepsilon_{ij} = \frac{1}{2}(u_{i,j} + u_{j,i})
    \label{equ:geometric_equation}
\end{equation}

\begin{equation}
    \sigma_{ij,j} + f_i = 0
    \label{equ:equilibrium_equation}
\end{equation}

\begin{equation}
    \begin{bmatrix}
        \sigma_{xx} \\
        \sigma_{yy} \\
        \sigma_{xy}
    \end{bmatrix}
    = \frac{E}{(1-2\mu)(1+\mu)}
    \begin{bmatrix}
        1-\mu & \mu & 0 \\
        \mu & 1-\mu & 0 \\
        0 & 0 & \frac{1-2\mu}{2}
    \end{bmatrix}
    \begin{bmatrix}
        \varepsilon_{xx} \\
        \varepsilon_{yy} \\
        \varepsilon_{xy}
    \end{bmatrix}
    \label{equ:physical_equation}
\end{equation}
\noindent
where all equations are expressed in Einstein notation, $i$ and $j$ subscripts with value of $1$ or $2$, $u_i$ denotes the displacement in $x_i$ axis, $\sigma_{ij}$ and $\varepsilon_{ij}$ denotes the stress and strain on surface $x_i$ in $x_j$ direction respectively, $f_i$ denotes the body force component in $x_i$ direction, $E$ denotes the Young's Modules, and $\mu$ denotes the Poisson's ratio.

\subsection{Dataset Generation}
\begin{figure}[!ht]
    \centering
    \includegraphics[width=0.7\linewidth]{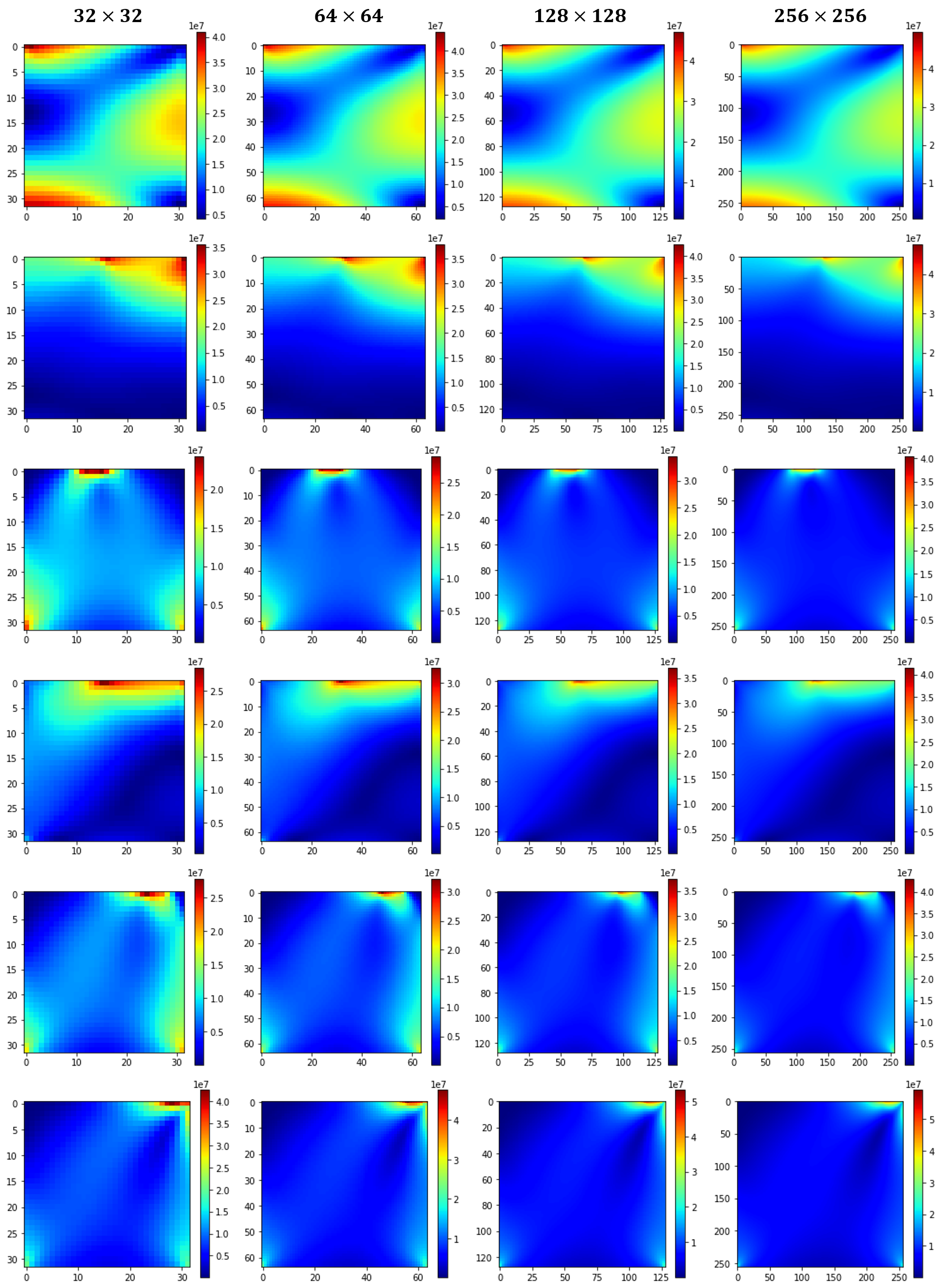}
    \caption{Dataset samples calculated by FEM. The mesh density from left to right is $32\times32$, $64\times64$, $128\times128$ and $256\times256$.}
    \label{fig:dataset_samples}
\end{figure}

We use the commercial FEA software ANSYS to generate the dataset required for SuperMeshingNet training and testing. The element node type is 4-node quadrilateral. Figure \ref{fig:dataset_samples} shows randomly selected samples from a total of 3240 samples. 

Each sample is composed of 4 images with different resolutions, which represent the von-mises stress field calculated by FEM with four mesh densities $24\times24, 48\times48, 128\times128$ and $256\times256$ from left to right. The neural network model inputs the stress field results of the $24\times24$ mesh density, and uses the stress field results of other mesh densities as ground truth. Among all samples, the von-mises stress varies from $0$ to $227.8$MPa with an average of $3.53$MPa.

Besides, in our experiments, we randomly rotate all samples during training to increase the diversity of the load and boundary conditions, improve the predictive ability of the model to deal with more input conditions.

\subsection{SuperMeshingNet Architecture}
\begin{figure*}[!ht]
    \centering
    \includegraphics[width=1\linewidth]{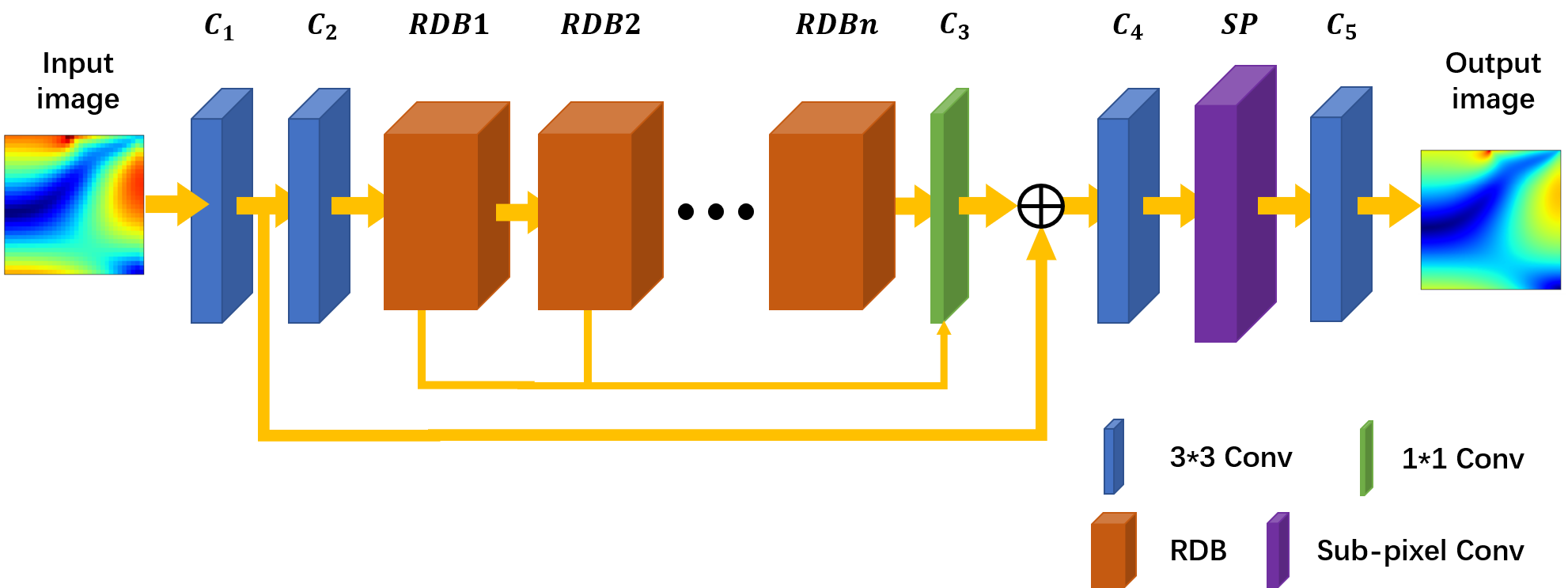}
    \caption{The architecture of SuperMeshingNet.}
    \label{fig:SuperMeshingNet}
\end{figure*}

As shown in Figure \ref{fig:SuperMeshingNet}, the proposed SuperMeshingNet increases the mesh density of the von-mises stress field. Model inputs low mesh-density stress field and outputs high mesh-density stress field, where each pixel represents the von-mises stress.

For SuperMeshingNet, we use the two convolutional layers $C_1$ and $C_2$ to extract the shallow feature of the input stress field image. Because of the similarity between the low and high mesh-density stress field, a global residual learning structure is introduced after $C_1$, so that the network can concentrate on the high-frequency information of the stress field to get better result.

After extracting the shallow feature, we utilize multiple RDBs for deep feature extraction. The deep feature in the low mesh-density stress field is extracted by cascading multiple RDBs. Then, concatenate the feature map extracted from each RDB and input into $C_3$. The feature extracted from each RDB can be further comprehensively used because of the concatenation. Convolutional layer $C_3$ has a filter with a kernel of $1\times1$, which is introduced to reduce the feature map dimension to reduce the network training difficulty.

Finally, we adopt a sub-pixel convolutional layer $SP$ for upsampling. Through the rearrangement of multiple feature maps, increase the scale of the feature map and obtain the stress field image with increased mesh density. 

In our approach, each convolutional layer has $64$ filters, and the kernel size of all filters is $3\times3$ except filters in $C_3$, which is $1\times1$.

\subsection{Loss Function and Metrics}
In our work, we select mean absolute error (MAE) as the loss function of our model for training. The predicted result of model $\hat{\sigma}$ and ground truth $\sigma$ are represented by image, according to different mesh-density scaling factor, the resolutions of the images are $64\times64$, $128\times128$ and $256\times256$ respectively. $\hat{\sigma}$ and $\sigma$ are reshaped into 1D vector before calculating MAE, $\hat{\sigma}$ is represented by $\hat{\sigma}=\lbrace\hat{\sigma}_1, \hat{\sigma}_2, \dots, \hat{\sigma}_n\rbrace$, and $\sigma$ is represented by $\sigma=\lbrace \sigma_1, \sigma_2, \dots, \sigma_n\rbrace$, MAE is: 

\begin{equation}
    MAE = \frac{1}{n}\sum_{i=1}^n \vert \hat{\sigma}_i - \sigma_i \vert
\end{equation}
\noindent
where $n$ is the amount of image pixels.

\section{Result and Discussion}
Our SuperMeshingNet is implemented by Pytorch, trained and tested on Nvidia GeForce GTX 2080Ti GPU. We carry out the following experiments and analyze the results to verify the feasibility of our method. First, we compare SuperMeshingNet with FEM in the case of mesh-density scaling factor is $8$, which prove the superiority of our method. Second, we test the effect of the number of different RDBs $D$ and the number of convolutional layers $C$ in them. Then, we analyze the effect of SuperMeshingNet on different mesh-density scaling factor. Finally, we test the effect of the model for the prediction of maximum stress.

\subsection{Accuracy and Performance}
\begin{figure}[!ht]
    \centering
    \includegraphics[width=0.7\linewidth]{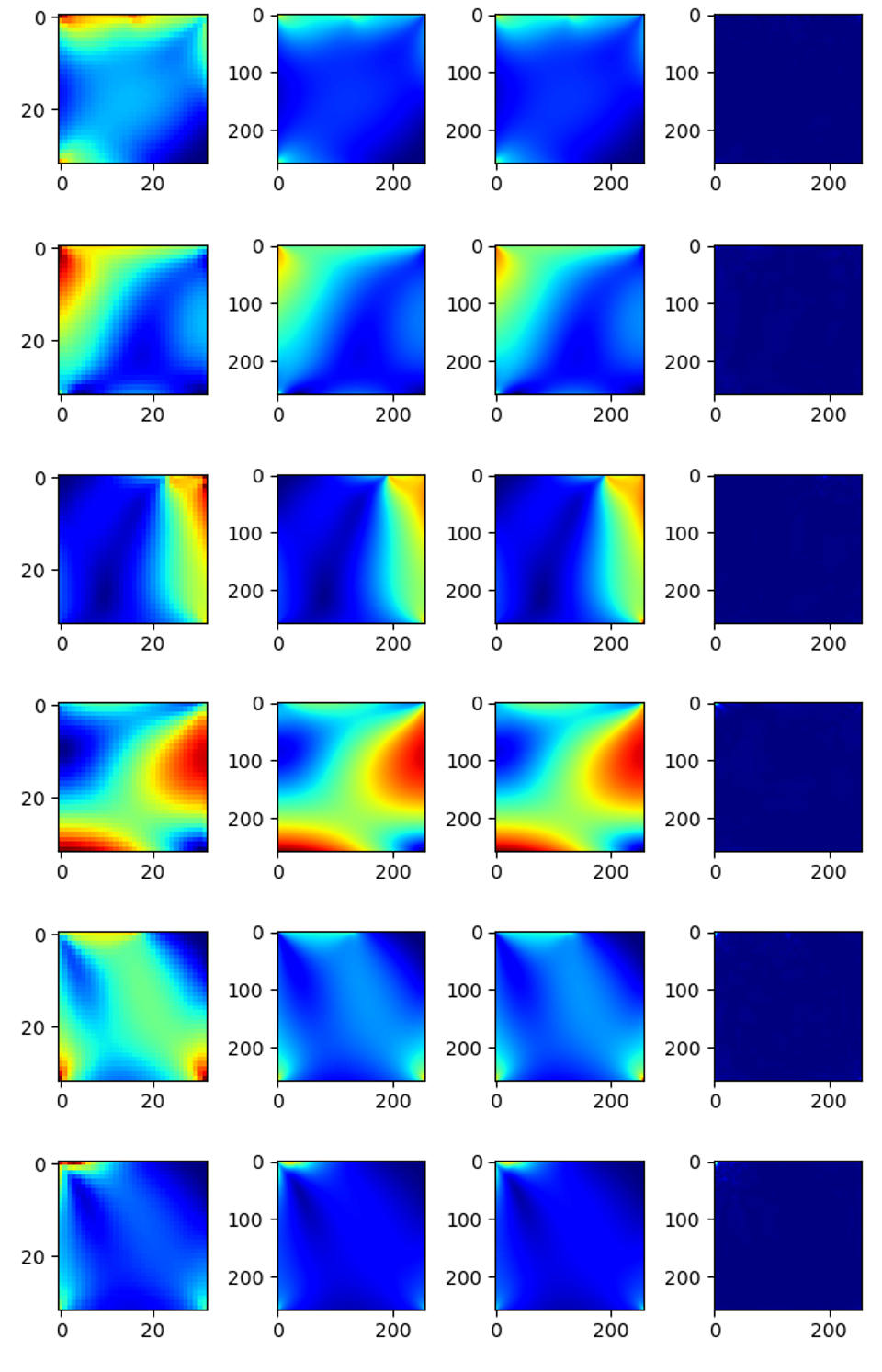}
    \caption{Several test results. the resolution of input image is $32\times32$, output image is $256\times256$.}
    \label{fig:test_results}
\end{figure}

In our experiments, we train and evaluate our SuperMeshingNet using the whole dataset. The training data size is 3000, and the test data size is 240. Figure \ref{fig:test_results} shows several test results under mesh-density scaling factor is $8$, the images from left to right are the input low mesh-density stress field, the input high mesh-density stress field, the generated high mesh-density stress field and the difference of the generated stress field and the ground truth separately. The average relative stress error (ARSE, defined as Equation (\ref{equ:ARSE})) between the generated high mesh-density stress field and the ground truth is only $0.54\%$, which we deem is within the acceptable range of stress field estimation. 

\begin{equation}
    ARSE = \frac{1}{nm} \sum_{i=1}^n \sum_{j=1}^m \vert \frac{\hat{\sigma}_i(j) - \sigma_i(j)}{\sigma_i(j)}\vert
    \label{equ:ARSE}
\end{equation}
\noindent
where $n$ denotes the amount of test cases, $m$ denotes the amount of image pixels. $\hat{\sigma}_i(j), \sigma_i(j)$ are the $jth$ stress value of the $ith$ test case of the generated high mesh-density stress field and the ground truth separately.

In order to verify the efficiency of our method, we test the computing efficiency between our method and FEM. Under mesh-density scaling factor is $8$, processing one sample by SuperMeshingNet in our method takes $0.715s$, and directly calculate high mesh-density stress field will cost $9.553s$. Our method is $13.36$ times faster than FEM at the same mesh density. All finite element computation is accomplished by commercial FEA software ANSYS (on Intel i5-7200U CPU).

\subsection{Effect of RDB Number and  Number of Convolutional Layers per RDB}
\begin{figure}[!ht]
    \centering
    \includegraphics[width=0.6\linewidth]{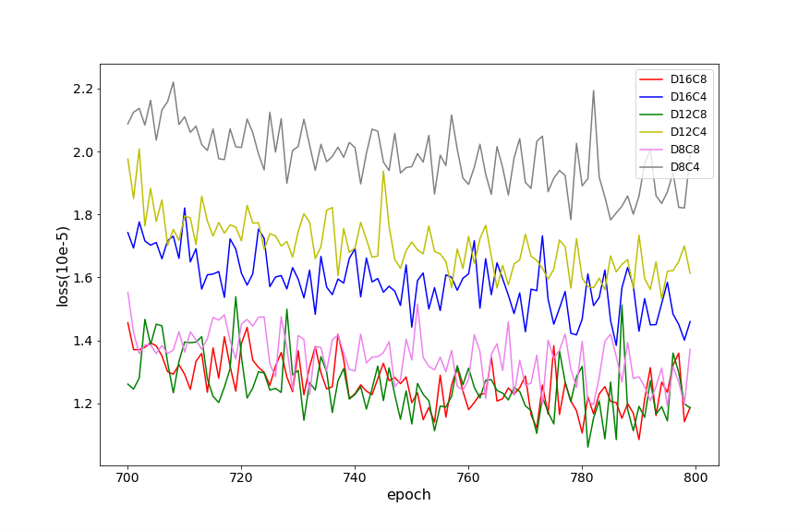}
    \caption{The MAE loss of different $D$ ($8, 12, 16$) and $C$ ($4, 8$) values in the last 100 epochs under $8$ mesh-density scaling factor.}
    \label{fig:D_C_effect}
\end{figure}

In order to obtain the best prediction performance, we test the effect of RDB number $D$ and number of convolutional layers per RDB $C$ through a series of experiments. We test a total of $6$ experiments, in which $D={8, 12, 16}$ and $C={4, 8}$ separately. Figure \ref{fig:D_C_effect} shows the MAE loss of these conditions under mesh-density scaling factor is $8$. Comparing the three pairs of MAE loss at the same $D$ value: $D16C8-D16C4, D12C8-D12C4, D8C8-D8C4$, the performance of SuperMeshingNet is significantly getting better with the increase of $C$. But the MAE loss of SuperMeshingNet at $D16C8, D12C8, D8C8$ are almost the same, which means $D$ has a slight effect on the model's performance. This result is different from the conclusion in \cite{SR_6}, which is larger $D$ and $C$ will lead to higher performance. We consider the reason to be SuperMeshingNet with smaller $D$ can perform well on this problem, so blindly increase $D$ will not further improve the performance. In our other experiments, we set $D=8$ and $C=8$ for satisfying performance and efficient training.

\subsection{Effect on Different Mesh-density Scaling Factor}
For the purpose of test the effect of our method under multiple mesh-density scaling factor $2, 4, 8$, we respectively calculated the ARSE between the generated high mesh-density stress field and the ground truth, as shown in Table \ref{tab:average relative stress error}. All the ARSE are in the acceptable range under different mesh-density scaling factor, which indicates our SuperMeshingNet has excellent prediction performance under multiple mesh-density scaling factor.

\begin{table}[!ht]
    \centering
    \renewcommand\arraystretch{1.5}
    \begin{tabular}{|c|c|c|c|}
        \hline
        Mesh density & $64\times64$ & $128\times128$ & $256\times256$ \\ \hline
        ARSE & 0.553\% & 0.547\% & 0.543\%  \\ \hline
    \end{tabular}
    \caption{The ARSE compared with the result calculated by FEM.}
    \label{tab:average relative stress error}
\end{table}

We also test the computing efficiency of our method, the computing time comparison of our method and FEM is shown in Table \ref{tab:computing efficiency}. We can find that the computing efficiency hardly decrease with the increase of mesh-density scaling factor, because the mapping from low mesh-density stress field to high mesh-density stress field has low time complexity relatively, and the computation of low mesh-density stress field is mainly time-consuming. Compared with the decrease of FEM, it shows that our method is more competitive in the computation of stress field with high mesh density.

\begin{table}[!ht]
    \centering
    \renewcommand\arraystretch{1.5}
    \begin{tabular}{|c|c|c|c|}
        \hline
        Mesh density    & $64\times64$ & $128\times128$ & $256\times256$ \\ \hline
        FEM+SuperMeshingNet      & 0.650 &  0.676  & 0.715   \\ \hline
        FEM                      & 1.132 &  3.176  & 9.553   \\ \hline
    \end{tabular}
    \caption{The computing time comparison of our method and FEM.}
    \label{tab:computing efficiency}
\end{table}

\subsection{Prediction of Maximum Stress}
\begin{figure}[!ht]
    \centering
    \includegraphics[width=1\linewidth]{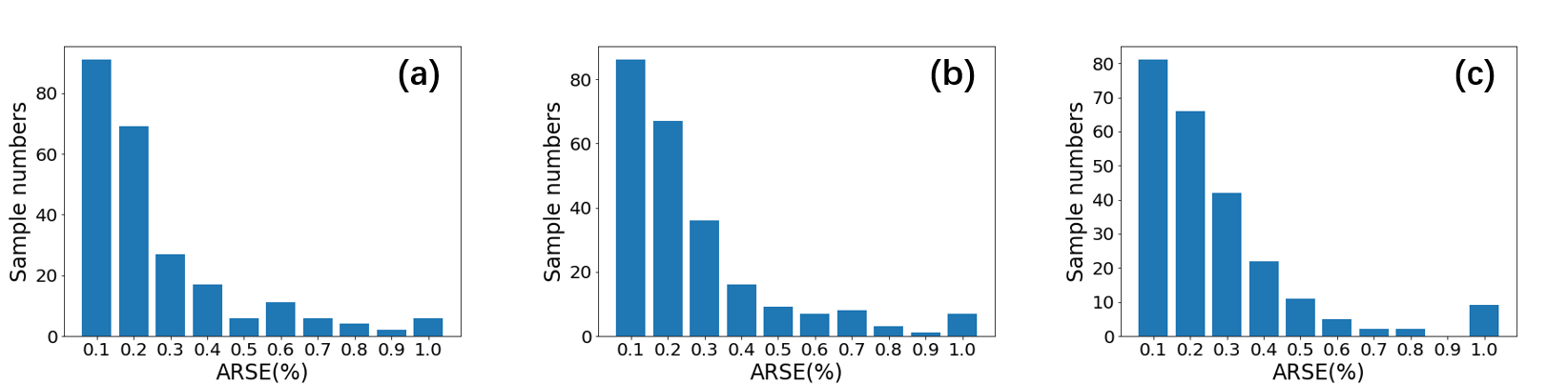}
    \caption{The distribution of the maximum stress difference between SuperMeshingNet predict result and FEM computing result. $(a)$: $64\times64$, $(b)$: $128\times128$, $(c)$: $256\times256$.}
    \label{fig:max_diff}
\end{figure}

For most plane strain problems, we mainly concerned about the maximum stress and its location to determine whether failure will occur at somewhere, so we test the prediction performance of maximum stress by our SuperMeshingNet. Figure \ref{fig:max_diff} shows the distribution of the maximum stress difference between generated high mesh-density stress field and ground truth on different mesh-density scaling factor. It's excellent that most of the maximum stress difference are under $0.4\%$, which proves the satisfying effect of SuperMeshingNet. However, we find that several cases' maximum difference is greater than $1.0\%$. We speculate the immature structure of SuperMeshingNet causes this, and we will continue to solve this problem in our future work. 

\section{Conclusions}
In this work, we present a novel method for boosting the mesh density in physical field prediction within 2D domain. In plane strain problems, we propose a new convolutional neural network called SuperMeshingNet mapping from low mesh-density stress field to high mesh-density stress field. In addition, we generate the training and testing dataset, which contains $3240$ samples of $32\times32, 64\times64, 128\times128, 256\times256$ mesh-density stress fields.

For the contradiction between spatial resolution and computing efficiency of FEA, we prove that our method can effectively predict physical field under various mesh densities through experiments. Compared to FEM, our method improve the computing efficiency $13.36$ times faster than FEM under the same spatial resolution, and the average relative stress error is only $0.54\%$, which means the computing accuracy is kept within the acceptable range of stress field prediction. 

In future work, we will focus on the unstructured mesh of the physical field relative to the image. Now we only take the structured mesh into consideration using CNN. Thus, it is more valuable that processing the more commonly used unstructured mesh by GNN. Moreover, integrating the inherent characteristics of the physical field (the strain continuity of the stress field, etc.) into neural network design can help improve prediction performance.

\bibliographystyle{unsrt}





\end{document}